\documentclass[9pt]{article}
\usepackage{mathrsfs}
\usepackage{amsthm}
\usepackage{amssymb}
\usepackage{amsmath}
\usepackage{graphicx}
\usepackage{color}
\usepackage{amsfonts}
\usepackage{float}
\usepackage{cite}
\usepackage [latin1]{inputenc}
\usepackage[text={152mm,234mm},left=38mm,vmarginratio=1:1]{geometry}
\newtheorem{theorem}{Theorem}[section]
\newtheorem{remark}{Remark}[section]
\newtheorem{lemma}[theorem]{Lemma}
\newtheorem{proposition}[theorem]{Proposition}

\numberwithin{equation}{section}
\normalsize

\begin{document}
\title{\textbf{Central limit theorems of occupation times of high-dimensional normalized binary contact path processes}}

\author{Xiaofeng Xue \thanks{\textbf{E-mail}: xfxue@bjtu.edu.cn \textbf{Address}: School of Mathematics and Statistics, Beijing Jiaotong University, Beijing 100044, China.}\\ Beijing Jiaotong University}

\date{}
\maketitle

\noindent {\bf Abstract:} The binary contact path process (BCPP) introduced in Griffeath (1983) describes the spread of an epidemic on a graph and is an auxiliary model in the study of improving upper bounds of the critical value of the contact process. In this paper, we are concerned with the central limit theorem of the occupation time of a normalized version of the BCPP (NBCPP) on a lattice. We show that the centred occupation time process of the NBCPP  converges in finite dimensional distributions to a Brownian motion when the dimension of the lattice and the infection rate of the model are sufficiently large and the initial state of the NBCPP is distributed with a particular invariant distribution.


\section{Introduction}\label{section one}
In this paper we are concerned with the normalized binary contact path process (NBCPP). For later use, we first introduce some notations. For $d\geq 1$, the $d$-dimensional lattice is denoted by $\mathbb{Z}^d$. For $x,y\in \mathbb{Z}^d$, we write $x\sim y$ when they are neighbors. The origin of $\mathbb{Z}^d$ is denoted by $O$. Now we recall the definition of the binary contact path process (BCPP) introduced in \cite{Grif1983} by Griffeath. The binary contact path process $\{\phi_t\}_{t\geq 0}$ on $\mathbb{Z}^d$ is a continuous-time Markov process with state space $\mathbb{X}=\{0,1,2,\ldots\}^{\mathbb{Z}^d}$ and evolves as follows. For any $x\in \mathbb{Z}^d$, $x\sim y$ and $t\geq 0$,
\[
\phi_t(x)\rightarrow
\begin{cases}
0 & \text{~at rate~}\frac{1}{1+2\lambda d},\\
\phi_t(x)+\phi_t(y) & \text{~at rate~}\frac{\lambda}{1+2\lambda d},
\end{cases}
\]
where $\lambda$ is a positive constant. As a result, the generator $\Omega$ of $\{\phi_t\}_{t\geq 0}$ is given by
\[
\Omega f(\phi)=\frac{1}{1+2\lambda d}\sum_{x\in \mathbb{Z}^d}\left(f(\phi^{x,-})-f(\phi)\right)+\frac{\lambda}{1+2\lambda d}\sum_{x\in \mathbb{Z}^d}\sum_{y\sim x}\left(f(\phi^{x,y})-f(\phi)\right)
\]
for any $f$ from $\mathbb{X}$ to $\mathbb{R}$ depending on finite coordinates and $x\in \mathbb{Z}^d$, where
\[
\phi^{x,-}(z)=
\begin{cases}
\phi(z) & \text{~if~}z\neq x,\\
0 & \text{~if~}z=x
\end{cases}
\]
and
\[
\phi^{x,y}(z)=
\begin{cases}
\phi(z) & \text{~if~}z\neq x,\\
\phi(x)+\phi(y) & \text{~if~}z=x
\end{cases}
\]
for all $z\in \mathbb{Z}^d$.

Intuitively, $\{\phi_t\}_{t\geq 0}$ describes the spread of an epidemic on $\mathbb{Z}^d$. The integer value a vertex taking is the seriousness of the illness on this vertex. A vertex taking value $0$ is healthy and one taking positive value is infected. An infected vertex becomes healthy at rate $\frac{1}{1+2\lambda d}$. A vertex $x$ is infected by a given neighbor $y$ at rate $\frac{\lambda}{1+2\lambda d}$. When the infection occurs, the seriousness of the illness on $x$ is added with that on $y$.

The BCPP $\{\phi_t\}_{t\geq 0}$ is introduced to improve upper bounds of critical values of high-dimensional contact processes (CP) according to the fact that the contact process $\{\xi_t\}_{t\geq 0}$ on $\mathbb{Z}^d$ can be equivalently defined as
\[
\xi_t(x)=
\begin{cases}
0 & \text{~if~} \phi_{(1+2\lambda d)t}(x)=0,\\
1 & \text{~if~} \phi_{(1+2\lambda d)t}(x)>0
\end{cases}
\]
for any $x\in \mathbb{Z}^d$. For a detailed survey of the contact process, see Chapter $6$ of \cite{Lig1985} and Part \uppercase\expandafter{\romannumeral2} of \cite{Lig1999}. The critical value $\lambda_c$ of the contact process is defined as
\[
\lambda_c=\sup\left\{\lambda:~P\left(\sum_{x\in \mathbb{Z}^d}\xi_t^{\lambda}(x)=0\text{~for some~}t>0\Bigg|\sum_{x}\xi^\lambda_0(x)=1\right)=1\right\}.
\]
According to the above coupling relationship between BCPP and CP, it is shown in \cite{Grif1983} that $\lambda_c(d)\leq \frac{1}{2d(2\gamma_d-1)}$ for $d\geq 3$, where $\gamma_d$  is the probability that the simple random walk on $\mathbb{Z}^d$ starting at $O$ never return to $O$ again.  Especially,
$\lambda_c(3)\leq 0.523$ as a corollary of the above upper bound. In \cite{Xue2018}, a modified version of BCPP is introduced and then a further
improved upper bound of $\lambda_c(d)$ is given for $d\geq 3$. It is shown in \cite{Xue2018} that $\lambda_c(d)\leq \frac{2-\gamma_d}{2d\gamma_d}$ and consequently $\lambda_c(3)\leq 0.340$.

The BCPP belongs to a family of continuous-time Markov processes called linear systems defined in Chapter 9 of \cite{Lig1985}, since there are a series of linear transformations $\{\mathcal{A}_k:~k\geq 1\}$ on $\mathbb{X}$ such that $\phi_t=\mathcal{A}_k\phi_{t-}$ for some $k\geq 1$ at each jump moment $t$.
As a result, for each $m\geq 1$, Kolmogorov-Chapman equations of
\[
\left\{\mathbb{E}\left(\prod_{i=1}^m\phi_t(x_i)\right):~x_1,\ldots,x_m\in \mathbb{Z}^d\right\}
\]
are given by a series of linear ordinary differential equations. For mathematical details, see Theorems 9.1.27 and 9.3.1 of \cite{Lig1985}.

For technical reasons which we will recall in Section \ref{section two}, it is convenient to investigate a normalized version $\{\eta_t\}_{t\geq 0}$ of the BCPP defined by
\[
\eta_t=\exp\left\{\frac{1-2d\lambda}{2d\lambda}t\right\}\phi_{\frac{1+2d\lambda}{2d\lambda}t}.
\]
The state space of the normalized binary contact path process (NBCPP) $\{\eta_t\}_{t\geq 0}$ is $\mathbb{Y}=[0, +\infty)^{\mathbb{Z}^d}$ and the generator $\mathcal{L}$ of $\{\eta_t\}_{t\geq 0}$ is given by
\begin{align*}
\mathcal{L}f(\eta)=&\frac{1}{2\lambda d}\sum_x\left(f(\eta^{x,-})-f(\eta)\right)+\frac{1}{2d}\sum_{x\in \mathbb{Z}^d}\sum_{y\sim x}\left(f(\eta^{x,y})-f(\eta)\right)\\
&+\left(\frac{1}{2\lambda d}-1\right)\sum_{x\in \mathbb{Z}^d}f^\prime_x(\eta)\eta(x),
\end{align*}
where $f_x^\prime$ is the partial derivative of $f$ with respect to the coordinate $\eta(x)$.

Inspired by counterpart results given in \cite{Cox1983} and \cite{Birkner2007} for voter models and branching random walks, in this paper we study central limit theorem of the occupation time
$\left\{\int_0^t\eta_u(O)du\right\}_{t\geq 0}$ of NBCPP. According to a calculation of the variance, it is natural to guess that the above central limit should be with three different forms in respective cases where $d\geq 5$, $d=4$ and $d=3$, as Theorems 1 of \cite{Cox1983} and 1.1 of \cite{Birkner2007}. However, in this paper we only give a rigorous result for part of the first case, i.e., CLT of the occupation time for sufficiently large $d$, according to the fact that we can not bound the fourth moment of $\eta_t(O)$ uniformly for $t\geq 0$ when the dimension $d$ is low. For our main result and mathematical details in the proof, see Sections \ref{section two} and \ref{section three}.

\section{Main results}\label{section two}
In this section we give our main result. For later use, we first introduce some nations and definitions. We denote by $\{S_n\}_{n\geq 0}$ the discrete-time simple random walk on $\mathbb{Z}^d$, i.e.,
\[
P\Big(S_{n+1}=y\Big|S_n=x\Big)=\frac{1}{2d}
\]
for any $n\geq 0$, $x\in \mathbb{Z}^d$ and $y\sim x$. We denote by $\{Y_t\}_{t\geq 0}$ the continuous-time simple random walk on $\mathbb{Z}^d$ with generator $L$ given by
\[
Lh(x)=\frac{1}{2d}\sum_{y\sim x}\left(h(y)-h(x)\right)
\]
for any bounded $h$ from $\mathbb{Z}^d$ to $\mathbb{R}$ and $x\in \mathbb{Z}^d$. As defined in Section \ref{section one}, we use $\gamma_d$ to denote the probability that $\{S_n\}_{n\geq 1}$ never return to $O$ again conditioned on $S_0=O$, i.e,
\[
\gamma_d=P\Big(S_n\neq O\text{~for all~}n\geq 1\Big|S_0=O\Big)=P\Big(Y_t\neq O\text{~for any~}t\geq 0\Big|Y_0=x\Big)
\]
for any $x\sim O$. For any $t\geq 0$, we use $p_t(\cdot, \cdot)$ to denote transition probabilities of $Y_t$, i.e.,
\[
p_t(x,y)=P\Big(Y_t=y\Big|Y_0=x\Big)
\]
for any $x,y\in \mathbb{Z}^d$. For any $x\in \mathbb{Z}^d$, we define
\[
\Phi(x)=P\Big(S_n=x\text{~for some~}n\geq 0\Big|S_0=O\Big).
\]
We use $\vec{1}$ to denote the configuration in $\mathbb{Y}$ where all vertices take value $1$. For any $\eta\in \mathbb{Y}$, we use $\mathbb{E}_\eta$ to denote the expectation operator of $\{\eta_t\}_{t\geq 0}$ conditioned on $\eta_0=\eta$. Furthermore, for any probability measure $\mu$ on $\mathbb{Y}$, we use $\mathbb{E}_\mu$ to denote the expectation operator of $\{\eta_t\}_{t\geq 0}$ conditioned on $\eta_0$ being distributed with $\mu$.

Now we recall some basic properties of NBCPP $\{\eta_t\}_{t\geq 0}$ proved in Chpater 9 of \cite{Lig1985} and Section 2 of \cite{Xue2021}.
\begin{proposition}\label{proposition 2.1}(Liggett, \cite{Lig1985})
For any $\eta\in \mathbb{Y}$, $t\geq 0$ and $x\in \mathbb{Z}^d$,
\[
\mathbb{E}_\eta\eta_t(x)=\sum_{y\in \mathbb{Z}^d}p_t(x, y)\mathbb{E}_\eta\eta_t(y).
\]
Furthermore, $\mathbb{E}_{\vec{1}}\eta_t(x)=1$ for any $x\in \mathbb{Z}^d$ and $t\geq 0$.
\end{proposition}

\begin{proposition}\label{proposition 2.2}(Liggett, \cite{Lig1985})
There exist a series of functions $\{q_t\}_{t\geq 0}$ from $(\mathbb{Z}^d)^2$ to $\mathbb{R}$ and a series of functions $\{\hat{q}_t\}_{t\geq 0}$ from $(\mathbb{Z}^d)^4$ to $\mathbb{R}$ such that
\[
\mathbb{E}_{\vec{1}}\left(\eta_t(O)\eta_t(x)\right)=\mathbb{E}_{\vec{1}}\left(\eta_t(y)\eta_t(x+y)\right)
=\sum_{z\in \mathbb{Z}^d}q_t(x,z)
\]
and
\[
\mathbb{E}_\eta\left(\eta_t(x)\eta_t(y)\right)=\sum_{z,w\in \mathbb{Z}^d}\hat{q}_t\left((x,y), (z,w)\right)\eta(z)\eta(w)
\]
for any $t\geq 0$, $x,y\in \mathbb{Z}^d$ and $\eta\in \mathbb{Y}$. Furthermore,
\[
q_t(O,w)=\sum_{(y,z):~y-z=w}\hat{q}_t\left((x,x), (y,z)\right)
\]
for any $x,z\in \mathbb{Z}^d$.
\end{proposition}
Proposition \ref{proposition 2.1} is an application of Theorem 9.1.27 of \cite{Lig1985} and Proposition \ref{proposition 2.2} is an application of Theorem 9.3.1 of \cite{Lig1985}. We write $q_t$ and $\hat{q}_t$ as $q_t^{\lambda, d}$ and $\hat{q}_t^{\lambda, d}$ respectively when we need to distinguish the dimension $d$ and infection rate $\lambda$.

\begin{proposition}\label{proposition 2.3}(Liggett, \cite{Lig1985})
When $d\geq 3$ and $\lambda>\frac{1}{2d(2\gamma_d-1)}$, we have following properties.

1)For any $x,y,z,w\in \mathbb{Z}^d$,
\[
\lim_{t\rightarrow+\infty}q_t^{\lambda, d}(x,y)=\lim_{t\rightarrow+\infty}\hat{q}_t^{\lambda, d}(x,y,z,w)=0.
\]

2) Conditioned on $\eta_0=\vec{1}$, $\eta_t$ converges weakly as $t\rightarrow+\infty$ to a probability measure $\nu_{\lambda, d}$ on $\mathbb{Y}$.

3) For any $x,y\in \mathbb{Z}^d$,
\[
\lim_{t\rightarrow+\infty}\mathbb{E}_{\vec{1}}\left(\left(\eta_t(x)\right)^2\right)=\sup_{t\geq 0}\mathbb{E}_{\vec{1}}\left(\left(\eta_t(x)\right)^2\right)
=\mathbb{E}_{\nu_{\lambda, d}}\left(\left(\eta_0(O)\right)^2\right)=1+\frac{1}{h_{\lambda, d}}
\]
and
\[
\lim_{t\rightarrow+\infty}{\rm Cov}_{\vec{1}}\left(\eta_t(x), \eta_t(y)\right)={\rm Cov}_{\nu_{\lambda, d}}\left(\eta_0(x), \eta_0(y)\right)=\frac{\Phi(x-y)}{h_{\lambda,d}},
\]
where $h_{\lambda, d}=\frac{2\lambda d(2\gamma_d-1)-1}{1+2d\lambda}$.
\end{proposition}

Proposition \ref{proposition 2.3} is an application of Corollary 2.8.20 and Theorem 9.3.17 of \cite{Lig1985}.

\begin{proposition}\label{proposition 2.4}(\cite{Xue2021})
There exist an integer $d_0\geq 5$ and a real number $\lambda_0>0$ satisfying following properties.

(1) If $d\geq d_0$ and $\lambda\geq \lambda_0$, then $\lambda>\frac{1}{2d(2\gamma_d-1)}$.

(2) For any $d\geq d_0$, $\lambda\geq \lambda_0$ and $x\in \mathbb{Z}^d$,
\[
\mathbb{E}_{\nu_{\lambda, d}}\left(\left(\eta_0(O)^4\right)\right)\leq\liminf_{t\rightarrow+\infty}\mathbb{E}^{\lambda, d}_{\vec{1}}\left(\left(\eta_t(x)\right)^4\right)<+\infty.
\]

(3) For any $d\geq d_0$ and $\lambda\geq \lambda_0$,
\[
\lim_{M\rightarrow+\infty}\sup_{(x,y,z)\in \left(\mathbb{Z}^d\right)^3:\atop \\\|x-y\|_1\wedge\|x-z\|_1\geq M}{\rm Cov}_{\nu_{\lambda, d}}\left(\left(\eta_0(x)\right)^2, \eta_0(y)\eta_0(z)\right)=0,
\]
where $\|\cdot\|_1$ is the $l_1$-norm on $\mathbb{R}^d$ and $a\wedge b=\min\{a,b\}$ for $a,b\in \mathbb{R}$.
\end{proposition}
Parts one and two of Proposition \ref{proposition 2.4} follow from Proposition 2.1 of \cite{Xue2021} and Fatou's lemma. Part three of Proposition \ref{proposition 2.4} follows from the analysis leading to Proposition 2.2 of \cite{Xue2021} and Fatou's lemma.

We omit proofs of Propositions \ref{proposition 2.1} to \ref{proposition 2.4} in this paper since they are only repeats of analyses given in Chapter 9 of \cite{Lig1985} and Sections 3, 4 of \cite{Xue2021}. Now we give our main result. For any $t\geq 0$ and $N\geq 1$, we define
\[
X_t^N=\frac{1}{\sqrt{N}}\int_0^{tN}\left(\eta_u(O)-1\right)du.
\]
We write $X_t^N$ as $X_{t,\lambda,d}^N$ when we need to distinguish $d$ and $\lambda$. We have the following theorem.

\begin{theorem}\label{theorem occupation main}
Let $d_0$ and $\lambda_0$ be defined as in Proposition \ref{proposition 2.4}. For any $d\geq d_0$, $\lambda>\lambda_0$, integer $m\geq 1$ and $t_1,t_2,\ldots,t_m\geq 0$, conditioned on the initial state $\eta_0$ of the NBCPP on $\mathbb{Z}^d$ being distributed with $\nu_{\lambda, d}$, $\left(X_{t_1, \lambda, d}^N, X_{t_2, \lambda, d}^N,\ldots, X_{t_m, \lambda, d}^N\right)$ converges weakly to $\sqrt{C_1(\lambda, d)}\left(B_{t_1}, B_{t_2},\ldots,B_{t_m}\right)$ as $N\rightarrow+\infty$, where $\{B_t\}_{t\geq 0}$ is a standard Brownian motion and
\[
C_1(\lambda, d)=\frac{2\int_0^{+\infty}\int_0^{+\infty}p_{r+\theta}(O, O)drd\theta}{h_{\lambda, d}\int_0^{+\infty}p_\theta(O, O)d\theta}.
\]
\end{theorem}

\begin{remark}\label{remark 2.6}
Theorem \ref{theorem occupation main} is consistent with a calculation of variance. According to Propositions \ref{proposition 2.1} to \ref{proposition 2.3} and the Markov property of $\{\eta_t\}_{t\geq 0}$,
\[
\lim_{N\rightarrow+\infty}{\rm Var}_{\nu_{\lambda, d}}(X_t^N)=\frac{2t}{h_{\lambda, d}}\int_0^{+\infty}\sum_x\Phi(x)p_r(O,x)dr.
\]
According to the strong Markov property of the simple random walk,
\[
\int_0^{+\infty}p_\theta(x,O)d\theta=\Phi(x)\int_0^{+\infty}p_\theta(O, O)d\theta
\]
and hence
\begin{align*}
\int_0^{+\infty}\sum_x\Phi(x)p_\theta(O,x)d\theta&=\frac{\int_0^{+\infty}\int_0^{+\infty}\sum_xp_r(O,x)p_\theta(x,O)drd\theta}{\int_0^{+\infty}p_\theta(O, O)d\theta}\\
&=\frac{\int_0^{+\infty}\int_0^{+\infty}p_{r+\theta}(O, O)drd\theta}{\int_0^{+\infty}p_\theta(O, O)d\theta}.
\end{align*}
\end{remark}

Theorem \ref{theorem occupation main} shows that the central limit theorem of the NBCPP is an analogue of that of voter models and branching random walks given in \cite{Cox1983} and \cite{Grif1983} when the dimension $d$ is sufficiently large. According to Proposition \ref{proposition 2.3}, calculations of variances as in Remark \ref{remark 2.6} imply that central limit theorems of occupation times of NBCPP on $\mathbb{Z}^3$, $\mathbb{Z}^4$ and $\mathbb{Z}^d$ for $d\geq 5$ should respectively be analogues of that of voter models and branching random walks in each cases. However, our current proof of Theorem \ref{theorem occupation main} relies heavily on the fact that $\sup_{t\geq 0}\mathbb{E}_{\vec{1}}^{\lambda, d}\left(\left(\eta_t(O)\right)^4\right)<+\infty$ when $d$ and $\lambda$ are sufficiently large, which we have not managed to prove yet for small $d$. That is why we currently only discuss the NBCPP on $\mathbb{Z}^d$ with $d$ sufficiently large.

The proof of Theorem \ref{theorem occupation main} is given in Section \ref{section three}, which follows the strategy introduced in \cite{Birkner2007} to prove central limit theorems of occupation times of branching random walks. The core idea of the strategy is to decompose $\int_0^t\left(\eta_u(O)-1\right)du$ as a martingale $M_t$ plus a remainder $R_t$ such that the quadratic variation process of $\frac{1}{\sqrt{N}}M_{tN}$ converges to $C_1 t$ in $L^2$ and $\frac{1}{\sqrt{N}}R_{tN}$ converges to $0$ in probability as $N\rightarrow+\infty$. To give the above decomposition, a resolvent function of the simple random walk $\{Y_t\}_{t\geq 0}$ on $\mathbb{Z}^d$ is utilized. For mathematical details, see Section \ref{section three}.

\section{Proof of Theorem \ref{theorem occupation main}}\label{section three}
In this section we prove Theorem \ref{theorem occupation main}. Throughout this section we assume that $d\geq d_0$ and $\lambda\geq \lambda_0$, where $d_0$ and $\lambda_0$ are defined as in Proposition \ref{proposition 2.4}. As we have introduced at the end of Section \ref{section two}, our proof follows the strategy introduced in \cite{Birkner2007}, where $X_t^N$ is decomposed as a martingale plus a remainder term. As $N\rightarrow+\infty$, the martingale converges weakly to a Brownian motion and the remainder converges to $0$ in probability. In detail, for any $\theta>0$, we define
\[
G_\theta(\eta)=\sum_{x\in \mathbb{Z}^d}g_\theta(x)\left(\eta(x)-1\right),
\]
where $g_\theta$ is the resolvent function of the simple random walk given by
\[
g_\theta(x)=\int_0^{+\infty}e^{-\theta u}p_u(O,x)du.
\]
According to the fact that $p_t(O, O)=O(t^{-d/2})$ as $t\rightarrow+\infty$, $\sum_{x} g_0^k(x)<+\infty$ when $k\geq 2$ and $d\geq 5$.
Let
\[
M_t^\theta=G_\theta(\eta_t)-G_\theta(\eta_0)-\int_0^t \mathcal{L}G_\theta(\eta_s)ds,
\]
then by Dynkin's martingale formula, $\{M_t^\theta\}_{t\geq 0}$ is a martingale with quadratic variation process $\{\langle M^\theta\rangle_t\}_{t\geq 0}$ given by
\[
\langle M^\theta\rangle_t=\int_0^t\left(\mathcal{L}\left(\left(G_\theta(\eta_s)\right)^2\right)-2G_\theta(\eta_s)\mathcal{L}G_\theta(\eta_s)\right)du.
\]
According to the definition of $\mathcal{L}$,
\begin{equation}\label{equ 3.1}
\langle M^\theta\rangle_t=\frac{1}{2\lambda d}\int_0^t\sum_{x\in \mathbb{Z}^d}\eta_s^2(x)\left(g_\theta^2(x)+\lambda\sum_{y\sim x}g_\theta^2(y)\right)ds.
\end{equation}

According to the fact that $\theta g_\theta(\cdot)-Lg_\theta(\cdot)=1_{O}(\cdot)$,
\[
\mathcal{L}G_\theta(\eta)=\theta G_\theta(\eta)-\left(\eta(O)-\theta\sum_xg_\theta(x)\right)=\theta G_\theta(\eta)-(\eta(O)-1).
\]
As a result,
\begin{equation}\label{equ 3.2}
X_t^N=\frac{1}{\sqrt{N}}M_{tN}^{1/N}+\frac{1}{\sqrt{N}}R_{tN}^{1/N},
\end{equation}
where
\[
R_t^\theta=-G_\theta(\eta_t)+G_\theta(\eta_0)+\int_0^t\theta G_\theta(\eta_s)ds.
\]
To prove Theorem \ref{theorem occupation main}, we need following lemmas.

\begin{lemma}\label{lemma 3.1}
For any $t\geq 0$, $\frac{1}{N}\langle M^{1/N}\rangle_{tN}$ converges to $C_1(\lambda, d)t$ in $L^2$ as $N\rightarrow+\infty$.
\end{lemma}

\begin{lemma}\label{lemma 3.2}
For any $t\geq 0$, $\frac{1}{\sqrt{N}}R_{tN}^{1/N}$ converges to $0$ in $L^2$ as $N\rightarrow+\infty$.
\end{lemma}

We first utilize Lemmas \ref{lemma 3.1} and \ref{lemma 3.2} to prove Theorem \ref{theorem occupation main}.

\proof[Proof of Theorem \ref{theorem occupation main}]
For any $t>0$, we claim that
\begin{equation}\label{equ uniformly integrable}
\lim_{N\rightarrow+\infty}\mathbb{E}\left(\sup_{0\leq s\leq tN}\frac{1}{N}\left(M^{1/N}_{s}-M^{1/N}_{s-}\right)^2\right)=0.
\end{equation}
We prove Equation \eqref{equ uniformly integrable} in Appendix \ref{subsection A.1}. By Equation \eqref{equ uniformly integrable}, Lemma \ref{lemma 3.1} and Theorem 1.4 of Chapter 7 of \cite{Ethier1986}, $\{\frac{1}{\sqrt{N}}M_{tN}^{1/N}\}_{t\geq 0}$ converges weakly to $\{\sqrt{C_1(\lambda, d)}B_t\}_{t\geq 0}$ as $N\rightarrow+\infty$. Consequently, Theorem \ref{theorem occupation main} follows from Equation \eqref{equ 3.2} and Lemma \ref{lemma 3.2}. 

\qed

At last, we give proofs of Lemmas \ref{lemma 3.1} and \ref{lemma 3.2}.

\proof[Proof of Lemma \ref{lemma 3.1}]
Conditioned on $Y_0=O$, the number of times $\{Y_t\}_{t\geq 0}$ visits $O$ follows a geometric distribution with parameter $\gamma_d$ and at each time $Y_t$ stays at $O$ for an exponential time with mean $1$. Therefore,
\[
\int_0^{+\infty}p_s(O,O)ds=\mathbb{E}_O\int_0^{+\infty}1_{\{X_s=O\}}ds=1\times \frac{1}{\gamma_d}=\frac{1}{\gamma_d}.
\]
Then by Equation \eqref{equ 3.1} and Proposition \ref{proposition 2.3},
\begin{align*}
\lim_{N\rightarrow+\infty}\mathbb{E}_{\nu_{\lambda, d}}\left(\frac{1}{N}\langle M^{1/N}\rangle_{tN}\right)&=\frac{t(1+2\lambda d)}{2\lambda d}\left(1+\frac{1}{h_{\lambda, d}}\right)\sum_{x}g_0^2(x)\\
&=\frac{2\gamma_d t}{h_{\lambda, d}}\int_0^{+\infty}\int_0^{+\infty}p_{\theta+r}(O, O)d\theta dr=C_1(\lambda, d)t.
\end{align*}
Hence, to complete the proof we only need to show that
\begin{equation}\label{equ 3.3}
\lim_{N\rightarrow+\infty}\frac{1}{N^2}{\rm Var}_{\nu_{\lambda,d}}\left(\langle M^{1/N}\rangle_{tN}\right)=0.
\end{equation}
By Proposition \ref{proposition 2.2} and the Markov property of $\{\eta_t\}_{t\geq 0}$, for $s<u$ and $x, w\in \mathbb{Z}^d$,
\[
{\rm Cov}_{\nu_{\lambda, d}}\left(\eta_s^2(x), \eta_u^2(w)\right)=\sum_{y,z}\hat{q}_{u-s}\left((w,w),(y,z)\right)
{\rm Cov}_{\nu_{\lambda, d}}\left(\eta_0^2(x), \eta_0(y)\eta_0(z)\right).
\]
Therefore, to prove Equation \eqref{equ 3.3} we only need to show that
\begin{equation}\label{equ 3.4}
\lim_{r\rightarrow+\infty}\sum_{x,w,y,z}V_0(x)V_0(w)\hat{q}_r\left((w,w), (y, z)\right)\left|{\rm Cov}_{\nu_{\lambda, d}}\left(\eta_0^2(x),\eta_0(y)\eta_0(z)\right)\right|=0,
\end{equation}
where $V_\theta(x)=g_\theta^2(x)+\sum_{y\sim x}g^2_\theta(y)$, which is decreasing in $\theta$. Since $\sum_xV_0(x)<+\infty$,
\[
\sum_{y,z}\hat{q}_r\left((w,w), (y,z)\right)=\sum_{y}q_r(O,y)=\mathbb{E}_{\vec{1}}\left(\left(\eta_r(O)\right)^2\right)
\leq 1+\frac{1}{h_{\lambda, d}}
\]
and
\[
\sup_{w,y,z}|{\rm Cov}_{\nu_{\lambda, d}}\left(\eta_0^2(w), \eta_0(y)\eta_0(z)\right)|<2\mathbb{E}_{\nu_{\lambda, d}}\left(\left(\eta_0(O)\right)^4\right)<+\infty
\]
according to Propositions \ref{proposition 2.2} to \ref{proposition 2.4} and Cauchy-Schwarz inequality, to prove \eqref{equ 3.4} we only need to show that
\begin{equation}\label{equ 3.5}
\lim_{r\rightarrow+\infty}\sum_{y,z}\hat{q}_r\left((w,w), (y,z)\right)\left|{\rm Cov}_{\nu_{\lambda, d}}\left(\eta_0^2(x), \eta_0(y)\eta_0(z)\right)\right|=0
\end{equation}
for any $x,w\in \mathbb{Z}^d$ according to the dominated convergence theorem. By Proposition \ref{proposition 2.4}, for any $\epsilon>0$, there exists $M>0$ such that $\left|{\rm Cov}_{\nu_{\lambda, d}}\left(\eta_0^2(x), \eta_0(y)\eta_0(z)\right)\right|<\epsilon$ when $\|x-y\|_1\wedge\|x-z\|_1\geq M$ and hence
\begin{align*}
&\sum_{(y,z):~\|x-y\|_1\wedge\|x-z\|_1\geq M}\hat{q}_r\left((w,w), (y,z)\right)\left|{\rm Cov}_{\nu_{\lambda, d}}\left(\eta_0^2(x), \eta_0(y)\eta_0(z)\right)\right|\\
&\leq \epsilon\sum_{(y,z):~\|x-y\|_1\wedge\|x-z\|_1\geq M}\hat{q}_r\left((w,w), (y,z)\right)\leq \epsilon(1+\frac{1}{h_{\lambda, d}}).
\end{align*}
We claim that
\begin{equation}\label{equ 3.6}
\lim_{t\rightarrow+\infty}\sum_{v}\hat{q}_t\left((w,w), (y,v)\right)=\lim_{t\rightarrow+\infty}\sum_{v}\hat{q}_t\left((w,w), (v,z)\right)=0
\end{equation}
for any $w,y,z\in \mathbb{Z}^d$. The proof of Equation \eqref{equ 3.6} is given in Appendix \ref{subsection A.2}. By Equation \eqref{equ 3.6}, Proposition \ref{proposition 2.4} and Cauchy-Schwarz inequality,
\[
\lim_{r\rightarrow+\infty}\sum_{(y,z):\|x-y\|_1\leq M\text{~or~}\|x-z\|_1\leq M}\hat{q}_r\left((w,w), (y,z)\right)\left|{\rm Cov}_{\nu_{\lambda, d}}\left(\eta_0^2(x), \eta_0(y)\eta_0(z)\right)\right|=0
\]
and hence
\[
\limsup_{r\rightarrow+\infty}\sum_{y,z}\hat{q}_r\left((w,w), (y,z)\right)\left|{\rm Cov}_{\nu_{\lambda, d}}\left(\eta_0^2(x), \eta_0(y)\eta_0(z)\right)\right|
\leq \epsilon(1+\frac{1}{h_{\lambda, d}}).
\]
Since $\epsilon$ is arbitrary, let $\epsilon\rightarrow 0$ and then Equation \eqref{equ 3.5} holds.

\qed

\proof[Proof of Lemma \ref{lemma 3.2}]

According to Propositions \ref{proposition 2.1} and \ref{proposition 2.3}, for any moment $c\geq 0$,
\begin{align*}
&\mathbb{E}_{\nu_{\lambda, d}}\left(\left(\frac{1}{\sqrt{N}}G_{1/N}(\eta_{c})\right)^2\right)={\rm Var}_{\nu_{\lambda, d}}\left(\frac{1}{\sqrt{N}}G_{1/N}(\eta_{c})\right)={\rm Var}_{\nu_{\lambda, d}}\left(\frac{1}{\sqrt{N}}G_{1/N}(\eta_0)\right)\\
&=\frac{1}{N}\sum_{x}\sum_{y}g_{1/N}(x)g_{1/N}(y)\frac{\Phi(y-x)}{h_{\lambda, d}}\\
&=\frac{\gamma_d}{Nh_{\lambda, d}}\int_0^{+\infty}\int_0^{+\infty}\int_0^{+\infty}e^{-(r+s)/N}\sum_x\sum_yp_r(0,x)p_s(x,y)p_u(y,0)drdsdu\\
&=\frac{\gamma_d}{Nh_{\lambda, d}}\int_0^{+\infty}\int_0^{+\infty}\int_0^{+\infty}e^{-(r+s)/N}p_{r+s+u}(O, O)drdsdu\\
&=\frac{\gamma_d}{Nh_{\lambda, d}}\int_0^{+\infty}p_{\theta}(O, O)\left(\int_0^{\theta}e^{-v/N}\left(\int_0^v1 du\right)dv\right)d\theta\\
&=\frac{N\gamma_d}{h_{\lambda, d}}\int_0^{+\infty}p_{\theta}(O, O)\left(1-e^{-\theta/N}\left(1+\frac{\theta}{N}\right)\right)d\theta.
\end{align*}
Since $p_\theta(O,O)=O(\theta^{-d/2})$ as $\theta\rightarrow+\infty$, $\lim_{x\rightarrow+\infty}1-e^{-x}(1+x)=1$ and $\lim_{x\rightarrow0}\frac{1-e^{-x}(1+x)}{x^2}=1/2$, there exist $0<M_1, C_4<+\infty$ such that $p_{\theta}(O, O)\leq C_4 \theta^{-d/2}$ when $\theta\geq M_1$, $1-e^{-x}(1+x)\leq C_4 x^2$ when $x\leq M_1$ and $1-e^{-x}(1+x)\leq C_4$ when $x>M_1$. Therefore, for $d\geq 5$,
\[
\frac{N\gamma_d}{h_{\lambda, d}}\int_0^{M_1}p_{\theta}(O, O)\left(1-e^{-\theta/N}\left(1+\frac{\theta}{N}\right)\right)d\theta
\leq \frac{C_4\gamma_d}{Nh_{\lambda, d}}\int_0^{M_1}\theta^2 d\theta=O(N^{-1}),
\]
\[
\frac{N\gamma_d}{h_{\lambda, d}}\int_{M_1}^{NM_1}p_{\theta}(O, O)\left(1-e^{-\theta/N}\left(1+\frac{\theta}{N}\right)\right)d\theta
\leq \frac{C_4^2\gamma_d}{Nh_{\lambda, d}}\int_{M_1}^{NM_1}\theta^{2-\frac{d}{2}}d\theta=O(N^{-1/2})
\]
and
\[
\frac{N\gamma_d}{h_{\lambda, d}}\int_{NM_1}^{+\infty}p_{\theta}(O, O)\left(1-e^{-\theta/N}\left(1+\frac{\theta}{N}\right)\right)d\theta
\leq \frac{C_4^2N\gamma_d}{h_{\lambda, d}}\int_{NM_1}^{+\infty}\theta^{-d/2}d\theta=O(N^{2-\frac{d}{2}}).
\]
In conclusion,
\begin{equation}\label{equ 3.8}
\lim_{N\rightarrow+\infty}\mathbb{E}_{\nu_{\lambda, d}}\left(\left(\frac{1}{\sqrt{N}}G_{1/N}(\eta_{0})\right)^2\right)
=\lim_{N\rightarrow+\infty}\mathbb{E}_{\nu_{\lambda, d}}\left(\left(\frac{1}{\sqrt{N}}G_{1/N}(\eta_{tN^2})\right)^2\right)=0.
\end{equation}
By Cauchy-Schwarz inequality,
\begin{align*}
\mathbb{E}_{\nu_{\lambda, d}}\left(\left(\frac{1}{\sqrt{N}}\int_0^{tN}\frac{1}{N}G_{1/N}(\eta_s)ds\right)^2\right)
&\leq\frac{tN}{N^3}\int_0^{tN}\left(\mathbb{E}_{\nu_{\lambda, d}}\left(\left(G_{1/N}(\eta_s)\right)^2\right)\right)ds\\
&=\frac{t^2N^2}{N^3}\mathbb{E}_{\nu_{\lambda, d}}\left(\left(G_{1/N}(\eta_0)\right)^2\right).
\end{align*}
Hence, according to Equation \eqref{equ 3.8},
\begin{equation}\label{equ 3.9}
\lim_{N\rightarrow+\infty}\mathbb{E}_{\nu_{\lambda, d}}\left(\left(\frac{1}{\sqrt{N}}\int_0^{tN}\frac{1}{N}G_{1/N}(\eta_s)ds\right)^2\right)=0.
\end{equation}
Lemma \ref{lemma 3.2} follows from Equations \eqref{equ 3.8} and \eqref{equ 3.9}.

\qed

\appendix{}
\section{Appendix}
\subsection{Proof of Equation \eqref{equ uniformly integrable}}\label{subsection A.1}

\proof[Proof of Equation \eqref{equ uniformly integrable}]

By Equation \eqref{equ 3.2} and the fact that $\{X_t^N\}_{t\geq 0}$ is continuous in $t$,
\begin{equation}\label{equ A.3}
\sup_{0\leq s\leq tN}\frac{1}{N}\left(M^{1/N}_{s}-M^{1/N}_{s-}\right)^2
\leq \sup_{0\leq s\leq tN}\frac{1}{N}\sup_{x}\max_{y\sim x}\{g^2_{1/N}(x)\eta^2_s(x), g^2_{1/N}(x)\eta^2_s(y)\}.
\end{equation}
For any $M>0$, let
\[
\tau_M=\inf\{s:~\frac{1}{N}\sup_{x}\max_{y\sim x}\{g^2_{1/N}(x)\eta^2_s(x), g^2_{1/N}(x)\eta^2_s(y)\}>M\},
\]
then
\[
\{\sup_{0\leq s\leq tN}\frac{1}{N}\sup_{x}\max_{y\sim x}\{g^2_{1/N}(x)\eta^2_s(x), g^2_{1/N}(x)\eta^2_s(y)\}>M\}=\{\tau_M\leq tN\}.
\]
Conditioned on $\{\tau_M\leq tN\}$, there exists $x_0\in \mathbb{Z}^d$ such that
\[
\max_{y\sim x_0}\{g^2_{1/N}(x_0)\eta^2_{\tau_M}(x_0), g^2_{1/N}(x_0)\eta^2_{\tau_ M}(y)\}>NM.
\]
If $\{\eta_s(y)\}_{y\sim x_0}$ and $\eta_s(x_0)$ do not jump to $0$ during $s\in [\tau_M, \tau_M+1]$, then
\[
\max_{y\sim x_0}\{g^2_{1/N}(x)\eta^2_{s}(x_0), g^2_{1/N}(x_0)\eta^2_{s}(y)\}\geq NMe^{\min\{(\frac{1}{2\lambda d}-1),0\}}
\]
for $s\in [\tau_M, \tau_{M+1}]$ and hence
\[
\int_0^{tN+1}\sum_{x}g_{1/N}^4(x)\left(\eta_s^4(x)+\sum_{y\sim x}\eta_s^4(y)\right)ds\geq N^2M^2e^{2\min\{(\frac{1}{2\lambda d}-1),0\}}.
\]
Since the state of a vertex jumps to $0$ at rate $1/(2\lambda d)$, according to the strong Markov property of $\{\eta_t\}_{t\geq 0}$,
\begin{align*}
&P\left(\int_0^{tN+1}\sum_{x}g_{1/N}^4(x)\left(\eta_s^4(x)+\sum_{y\sim x}\eta_s^4(y)\right)ds\geq N^2M^2e^{2\min\{(\frac{1}{2\lambda d}-1),0\}}\right)\\
&\geq P(\tau_M\leq tN)e^{-\frac{2d+1}{2\lambda d}}.
\end{align*}
Then, by Markov's inequality, Proposition \ref{proposition 2.4} and the fact that $\sum_x g_{1/N}^4(x)\leq \sum_x g_0^4(x)<+\infty$,
\begin{align}\label{equ A.4}
&P\left(\sup_{0\leq s\leq tN}\frac{1}{N}\sup_{x}\max_{y\sim x}\{g^2_{1/N}(x)\eta^2_s(x), g^2_{1/N}(x)\eta^2_s(y)\}>M\right) \notag\\
&\leq \frac{e^{\frac{2d+1}{2\lambda d}}}{N^2M^2e^{2\min\{(\frac{1}{2\lambda d}-1),0\}}}\int_0^{tN+1}\sum_xg_0^4(x)\left(\mathbb{E}_{\lambda, d}\left(\eta^4(O)\right)(2d+1)\right)ds \leq \frac{C_3}{NM^2},
\end{align}
where $C_3<+\infty$ is independent of $N$ and $M$. According to Fubini theorem, for any positive random variable $V$ and any $c>0$,
\[
\mathbb{E}\left(V1_{\{V>c\}}\right)=cP(V>c)+\int_c^{+\infty}P(V\geq u)du.
\]
Hence, for any $\epsilon>0$,
\[
\mathbb{E}\left(\sup_{0\leq s\leq tN}\frac{1}{N}\left(M^{1/N}_{s}-M^{1/N}_{s-}\right)^2\right)
\leq \epsilon+\epsilon\frac{C_3}{N\epsilon^2}+\int_{\epsilon}^{+\infty}\frac{C_3}{Nu^2}du
=\epsilon+\frac{2C_3}{N\epsilon}
\]
according to Equations \eqref{equ A.3} and \eqref{equ A.4}. As a result,
\[
\limsup_{N\rightarrow+\infty}\mathbb{E}\left(\sup_{0\leq s\leq tN}\frac{1}{N}\left(M^{1/N}_{s}-M^{1/N}_{s-}\right)^2\right)\leq \epsilon.
\]
Since $\epsilon$ is arbitrary, let $\epsilon\rightarrow 0$ and the proof is complete.

\qed

\subsection{Proof of Equation \eqref{equ 3.6}}\label{subsection A.2}

\proof[Proof of Equation \eqref{equ 3.6}]

By Theorem 9.3.1 of \cite{Lig1985},
\[
\hat{q}_t\left((w,w), (y,v)\right)=e^{-2t}\sum_{n=0}^{+\infty}\frac{t^n}{n!}H^n\left((w,w), (y,v)\right),
\]
where $H$ is a $(\mathbb{Z}^d)^2\times (\mathbb{Z}^d)^2$ matrix given by
\[
H((x,y), (u,v))=
\begin{cases}
\frac{1}{2d} &\text{~if~}x\neq y, u\sim x\text{~and~}v=y,\\
\frac{1}{2d} &\text{~if~}x\neq y, u=x \text{~and~}v\sim y,\\
\frac{1}{2d\lambda} & \text{~if~}x=x \text{~and~}(u,v)=(x,x),\\
\frac{1}{2d} & \text{~if~}x=x, u\sim x\text{~and~}v=u,\\
\frac{1}{2d} & \text{~if~}x=x, u=x\text{~and~}v\sim x,\\
\frac{1}{2d} & \text{~if~}x=x, u\sim x\text{~and~}v=x,\\
0 & \text{~else}
\end{cases}
\]
and
\begin{align*}
H^n\left((w,w), (y,v)\right)
=\sum_{\{(u_i, v_i)\}_{0\leq i\leq n}:\atop ~(u_0, v_0)=(w,w), (u_n, v_n)=(y,v)}\prod_{i=0}^{n-1}H\left((u_i, v_i), (u_{i+1}, v_{i+1})\right).
\end{align*}
As a result,
\[
\sum_{v\in \mathbb{Z}^d}\hat{q}_t\left((w,w), (y,v)\right)
=e^{-2t}\sum_{n=0}^{+\infty}\frac{(2t)^n}{n!}
\mathbb{E}_{(w,w)}\left(\prod_{i=0}^{n-1}\Theta\left(\beta_{i},\beta_{i+1}\right)1_{\{\beta_n(1)=y\}}\right),
\]
where $\{\beta_n=(\beta_n(1), \beta_n(2)):~n\geq 0\}$ is a random walk on $(\mathbb{Z}^d)^2$ such that
\[
P\Big(\beta_{n+1}=(u,v)\Big|\beta_n=(x,y)\Big)=
\begin{cases}
\frac{1}{2d} &\text{~if~}x\neq y, u\sim x\text{~and~}v=y,\\
\frac{1}{2d} &\text{~if~}x\neq y, u=x\text{~and~}v\sim y, \\
\frac{1}{6d+1} &\text{~if~}x=y\text{~and~}(u,v)=(x,x),\\
\frac{1}{6d+1} &\text{~if~}x=y, u\sim x\text{~and~}v=u,\\
\frac{1}{6d+1} &\text{~if~}x=y, u\sim x\text{~and~}v=x,\\
\frac{1}{6d+1} &\text{~if~}x=y, u=x{~and~}v\sim x,\\
0 &\text{~else}
\end{cases}
\]
and $\Theta$ is a function from $(\mathbb{Z}^d)^2$ to $\mathbb{R}$ such that
\[
\Theta\left((x,y), (u,v)\right)=\frac{H\left((x,y), (u,v)\right){\rm deg}(x,y)}{2},
\]
where
\[
{\rm deg}(x,y)=
\begin{cases}
4d & \text{~if~}x\neq y,\\
6d+1 & \text{~if~}x=y.
\end{cases}
\]
Then for each $i\geq 0$, $\Theta(\beta_i, \beta_{i+1})\geq 1$ and $\Theta(\beta_i, \beta_{i+1})=1$ when and only when $\beta_i(1)=\beta_i(2)$. Hence, by H\"{o}lder's inequality,
\begin{align}\label{equ A.1}
&\sum_{v\in \mathbb{Z}^d}\hat{q}_t\left((w,w), (y,v)\right) \notag\\
&\leq \left(\mathbb{E}_{(w,w)}\left(\prod_{i=0}^{+\infty}\Theta^{1+\epsilon}\left(\beta_{i},\beta_{i+1}\right)\right)\right)^{\frac{1}{1+\epsilon}}
e^{-2t}\sum_{n=0}^{+\infty}\frac{(2t)^n}{n!}\left(P\left(\beta_n(1)=y\right)\right)^{\epsilon/(1+\epsilon)}
\end{align}
for any $\epsilon>0$. Since $\{\beta_n(1)-\beta_n(2)\}_{n\geq 0}$ is a lazy version of the simple random walk on $\mathbb{Z}^d$ and $H(\beta_i, \beta_{i+1})>1$ only when $\beta_i(1)=\beta_i(2)$, by the strong Markov property of $\{\beta_n\}_{n\geq 0}$,
\[
\mathbb{E}_{(w,w)}\left(\prod_{i=0}^{+\infty}\Theta^{1+\epsilon}\left(\beta_{i},\beta_{i+1}\right)\right)
=\frac{4d}{6 d+1}\left(\frac{6d+1}{4d}\right)^{1+\epsilon}\gamma_d+\sum_{k=1}^{+\infty}\left(\alpha_\epsilon(d,\lambda)\right)^k\gamma_d,
\]
where
\[
\alpha_\epsilon(d,\lambda)=\frac{1}{6d+1}\left(\frac{6d+1}{2\cdot 2d\lambda}\right)^{1+\epsilon}
+\frac{2d}{6d+1}\left(\frac{6d+1}{2\cdot 2d}\right)^{1+\epsilon}
+\frac{4d}{6d+1}\left(\frac{6d+1}{2\cdot 2d}\right)^{1+\epsilon}(1-\gamma_d).
\]
Since $\lambda>\frac{1}{2d(2\gamma_d-1)}$,
\[
\lim_{\epsilon\rightarrow 0}\alpha_\epsilon(d,\lambda)=\frac{1}{4d\lambda}+1/2+1-\gamma_d<1
\]
and hence there exists $\epsilon_0>0$ such that $\mathbb{E}_{(w,w)}\left(\prod_{i=0}^{+\infty}\Theta^{1+\epsilon_0}\left(\beta_{i},\beta_{i+1}\right)\right)<+\infty$.  By Equation \eqref{equ A.1}, to complete the proof we only need to show that
\begin{equation}\label{equ A.2}
\lim_{t\rightarrow+\infty}e^{-2t}\sum_{n=0}^{+\infty}\frac{(2t)^n}{n!}\left(P\left(\beta_n(1)=y\right)\right)^{\epsilon_0/(1+\epsilon_0)}=0.
\end{equation}
Since $d\geq 3$ and $\{\beta_n(1)\}_{n\geq 1}$ is a lazy version of the simple random walk on $\mathbb{Z}^d$,
\[
\lim_{n\rightarrow+\infty}P\left(\beta_n(1)=y\right)=0.
\]
Hence, for any $\epsilon>0$, there exists $M\geq 1$ such that $\left(P\left(\beta_n(1)=y\right)\right)^{\epsilon_0/(1+\epsilon_0)}\leq \epsilon$ when $n\geq M$ and then
\[
e^{-2t}\sum_{n=M}^{+\infty}\frac{(2t)^n}{n!}\left(P\left(\beta_n(1)=y\right)\right)^{\epsilon_0/(1+\epsilon_0)}\leq \epsilon e^{-2t} \sum_{n=0}^{+\infty}\frac{(2t)^n}{n!}=\epsilon.
\]
Since $\lim_{t\rightarrow+\infty}\sum_{n=0}^{M-1}e^{-2t}\frac{(2t)^n}{n!}=0$, we have
\[
\limsup_{t\rightarrow+\infty}e^{-2t}\sum_{n=0}^{+\infty}\frac{(2t)^n}{n!}\left(P\left(\beta_n(1)=y\right)\right)^{\epsilon_0/(1+\epsilon_0)}\leq \epsilon.
\]
Since $\epsilon$ is arbitrary, let $\epsilon\rightarrow0$ and then Equation \eqref{equ A.2} holds.

\qed

\quad

\textbf{Acknowledgments.}
The author is grateful to the financial
support from Beijing Jiaotong University with grant number 2022JBMC039.

{}

\begin{thebibliography}{}
\bibitem{Birkner2007}Birkner, M. and Z\"{a}hle, I. (2007). A functional CLT for the occupation time
of a state-dependent branching random walk. \emph{The Annals of Probability} \textbf{35}, 2063-2090.
\bibitem{Cox1983}Cox, J. T. and Griffeath, D. (1983).
Occupation time limit theorems for the voter model. \emph{The Annals of Probability} \textbf{11}, 876-893.
\bibitem{Ethier1986}Ethier, N. and Kurtz, T. (1986). \emph{Markov Processes: Characterization and Convergence.} John Wiley and Sons, Hoboken, NJ, USA.
\bibitem{Grif1983}Griffeath, D. (1983). The Binary Contact Path Process. \emph{The Annals of Probability} \textbf{11}, 692-705.
\bibitem{Lig1985}Liggett, T. M. (1985). \emph{Interacting Particle Systems.} Springer, New York.
\bibitem{Lig1999}Liggett, T. M. (1999). \emph{Stochastic interacting systems: contact, voter and exclusion processes.}
Springer, New York.
\bibitem{Xue2018}Xue, XF. (2018). An improved upper bound for critical value of
the contact process on $\mathbb{Z}^d$ with $d\geq 3$. \emph{Electronic Communications in Probability} \textbf{23}, 77, 1-11.
\bibitem{Xue2021}Xue, XF. and Zhao, LJ. (2021). Non-equilibrium fluctuations of the weakly asymmetric normalized binary contact path process. \emph{Stochastic Processes and their Applications} \textbf{135}, 227-253.
\end{thebibliography}
\end{document}